# Lie Subalgebras and Invariant Solutions to the Equation of Fluid Flows in Toroidal Field

Lang Xia

Email: langxia.org@gmail.com

**Abstract:** Partial differential equations (PDEs), particularly coupled PDE systems, are difficult to solve. By far, no systematical approach is available to guarantee analytical solutions to an arbitrary PDE. The Lie symmetry method that is the generalization of Galois's approach for algebraic equations, however, has the potential to analyze and solve a class of differential equations systematically. In the present report, by using the Stokes-Helmholtz decomposition theorem the 3-dimensional Navier-Stokes equation (NSE) is uncoupled and transformed into a scalar equation for the velocity potential when the flow field is toroidal. The dynamics of the velocity potential is independent of the vector potential. The reduction and invariant solutions to the equation are analyzed by the Lie symmetry method subsequently. The Lie subalgebras for the equation are discussed and the corresponding invariant solutions are also presented.

Keywords: Navier-Stokes equation, Geometric analysis, PDE, Lie groups, Lie symmetries

## Introduction

The Navier-Stokes equation (NSE) in the classic mechanics is probably the most important PDE. Due to its complexity, neither the uniqueness nor existence of global solutions has been proved until now [1]. Although very few analytical solutions are well-known [2], researches in the fluid mechanics whose phenomena are governed by the NSE still advance fast with the help of numerical computations [3-6].

However, analytical solutions to the NSE are still precious. They provide a gauge for the verification and validation of all kind of numerical schemes and algorithms. Therefore, many semi-analytical methods and alternative analysis are developed to study the property of possible analytical solutions [7, 8].

Enlightened by the success of Galois groups in solving algebraic equations, Sophus Lie proposed continuous groups for solving differential equations [9, 10]. This method helped Lie to discover the symmetries hidden in the differential equations. Meanwhile, it unified many kinds of *ad hoc* methods of reducing differential equations. Although Lie groups have richer meanings in the fields of differential geometry and mathematical physics, It is only after 1970 they become a powerful tool to simplify and solve differential equations [11, 12]. That is because the Lie group method requires an enormous amount of algebraic calculations, and sometimes solving the linear determining equations are even harder than solving the original equations. Thanks to the developments of modern computer systems, those determining equations can be easily solved with the help of CAS software or packages right now. Therefore, within the past 30 years, the symmetric groups, Lie algebras and the corresponding invariant solutions to many famous equations have been gradually documented [2, 13].

In this paper, we focus on the simplification and Lie group analysis of the NSE. Firstly, we derive a relatively simple scalar equation using the Stokes-Helmholtz decomposition when the flow field is toroidal. Then the Lie group analysis is carried on to the equation with the assistance of the machinery of jet bundles. Finally, we present some solved Lie algebras and invariant solutions to the equation for special cases.

## Simplification of the NSE

Let consider the 3-D non-dimensional NSE without external forces in a simply connected domain that is written in the form of

$$\frac{\partial \mathbf{u}}{\partial t} + \mathbf{u} \cdot \nabla \mathbf{u} = -\nabla p + \frac{1}{\text{Re}} \nabla^2 \mathbf{u} \qquad (1)$$
$$\nabla \cdot \mathbf{u} = 0$$

where $p$ is the mechanical pressure, Re the Reynolds number and the velocity $\mathbf{u} \in C^\infty(\mathbb{R}^3)$. There is more than one approach to do the kinematic decomposition of flow fields (vector fields) [8], but the Stokes-Helmholtz might be the most widely used one. The purpose of this section is to derive a dynamical equation governing the flow field. Recall the Stokes-Helmholtz decomposition theorem, we write the velocity in the form of

$$\mathbf{u} = \nabla \phi + \nabla \times \boldsymbol{\psi} \qquad (2)$$

where $\phi$ the potential component is irrotational, and $\boldsymbol{\psi}$ the rotational component is solenoidal. Here we also require $\phi$ being an analytical function so that the decomposition is valid in any simply connected unbounded domain. It is well known that the decomposition is not unique [14], because adding any harmonic function to the rotational component will result in a new decomposition to the velocity field. This issue can be resolved physically by imposing a boundary condition such that the derivative of the velocity is bounded, more specifically, let $\mathbf{n}$ be the normal unit vector of the boundary, then the Stokes-Helmholtz decomposition is unique if and only if $\mathbf{n} \cdot (\nabla \times \boldsymbol{\psi}) = 0$ (in this case, we may call it Helmholtz–Hodge decomposition) [15]. Furthermore, by defining an equivalent class for the potential or rotational component $\nabla\phi/\sim$ mathematically, that is to say, $\nabla\phi \equiv \nabla\phi^1 \pmod{\nabla\Phi}$, here $\Phi$ is an arbitrary harmonic function we can also avoid the ambiguity. We shall see the discussion in our future paper. However, here we focus on the problem freed of boundary conditions and treat the decomposition "roughly" unique. To deal with the vector potential, let's recall the poloidal-toroidal decomposition (also known as Mie representation) of the solenoidal part in the spherical coordinate $(r,\theta,\varphi)$ [16]

$$\nabla \times \boldsymbol{\psi} = \mathbf{T} + \mathbf{P} = \nabla \times \left(\psi \frac{\mathbf{r}}{r}\right) + \nabla \times \nabla \times \left(\sigma \frac{\mathbf{r}}{r}\right) \tag{3}$$

In Cartesian coordinate, the above equation has a slightly complicated expression with mean fields being added [17]. We shall think the case that the poloidal field vanishes. Assume that $\mathbf{r} = (\mathbf{i}+\mathbf{j}+\mathbf{k})r$, and express the toroidal field in Cartesian coordinate locally as

$$\mathbf{T} = \nabla \times (\psi(\mathbf{i}+\mathbf{j}+\mathbf{k})) \tag{4}$$

Therefore, we can assume the vector potential $\boldsymbol{\psi} = (\psi,\psi,\psi)$, and then the three components of the velocity have the following forms

$$\begin{aligned} u_1 &= \frac{\partial \phi}{\partial x} + \frac{\partial \psi}{\partial y} - \frac{\partial \psi}{\partial z} \\ u_2 &= \frac{\partial \phi}{\partial y} + \frac{\partial \psi}{\partial z} - \frac{\partial \psi}{\partial x} \\ u_3 &= \frac{\partial \phi}{\partial z} + \frac{\partial \psi}{\partial x} - \frac{\partial \psi}{\partial y} \end{aligned} \tag{5}$$

The Eq.(5) indicate the vector potential is a toroidal field. From the above analysis, we have the following result

**Theorem:** If the NSE has the solution in the form of $\mathbf{u} = \nabla\phi + \nabla\times\boldsymbol{\psi}$, where the vector potential is a toroidal field, then the harmonic function satisfies the following equation

$$\frac{\partial \phi}{\partial t} + \frac{1}{3}\left(\frac{\partial \phi}{\partial x} + \frac{\partial \phi}{\partial y} + \frac{\partial \phi}{\partial z}\right)^2 - \frac{1}{\mathrm{Re}}\left(\frac{\partial^2 \phi}{\partial x \partial y} + \frac{\partial^2 \phi}{\partial y \partial z} + \frac{\partial^2 \phi}{\partial z \partial x}\right) + p = 0 \tag{6}$$

Proof: The summation of the velocity components equations gives rise to

$$u_1 + u_2 + u_3 = \frac{\partial \phi}{\partial x} + \frac{\partial \phi}{\partial y} + \frac{\partial \phi}{\partial z} \tag{7}$$

Note that the NSE also has the following component forms

$$\begin{aligned}\frac{\partial u_1}{\partial t} + \mathbf{u} \cdot \nabla u_1 &= -\frac{\partial p}{\partial x} + \frac{1}{Re}\nabla^2 u_1 \\ \frac{\partial u_2}{\partial t} + \mathbf{u} \cdot \nabla u_2 &= -\frac{\partial p}{\partial y} + \frac{1}{Re}\nabla^2 u_2 \\ \frac{\partial u_3}{\partial t} + \mathbf{u} \cdot \nabla u_3 &= -\frac{\partial p}{\partial z} + \frac{1}{Re}\nabla^2 u_3\end{aligned} \tag{8}$$

After the summation of the above equations, substitute Eq.(7) into it, we have

$$\frac{\partial}{\partial t}(\frac{\partial \phi}{\partial x} + \frac{\partial \phi}{\partial y} + \frac{\partial \phi}{\partial z}) + \mathbf{u} \cdot \nabla(\frac{\partial \phi}{\partial x} + \frac{\partial \phi}{\partial y} + \frac{\partial \phi}{\partial z}) = -(\frac{\partial p}{\partial x} + \frac{\partial p}{\partial y} + \frac{\partial p}{\partial z}) + \frac{1}{Re}\nabla^2(\frac{\partial \phi}{\partial x} + \frac{\partial \phi}{\partial y} + \frac{\partial \phi}{\partial z}) \tag{9}$$

For a scalar function $f$, we know that $\nabla \cdot (f\mathbf{u}) = f\nabla \cdot \mathbf{u} + \mathbf{u}\nabla f$. Apply the incompressible condition $\nabla \cdot \mathbf{u} = 0$, we get

$$\mathbf{u} \cdot \nabla(\frac{\partial \phi}{\partial x} + \frac{\partial \phi}{\partial y} + \frac{\partial \phi}{\partial z}) = \nabla \cdot (\mathbf{u}(\frac{\partial \phi}{\partial x} + \frac{\partial \phi}{\partial y} + \frac{\partial \phi}{\partial z})) \tag{10}$$

Let

$$\begin{aligned}\nabla \cdot \boldsymbol{\varphi} &= \frac{\partial \phi}{\partial x} + \frac{\partial \phi}{\partial y} + \frac{\partial \phi}{\partial z} \\ \nabla \cdot \mathbf{p} &= \frac{\partial p}{\partial x} + \frac{\partial p}{\partial y} + \frac{\partial p}{\partial z}\end{aligned} \tag{11}$$

where $\boldsymbol{\varphi} := \phi(x, y, z, t)(\mathbf{i}+\mathbf{j}+\mathbf{k})$ and $\mathbf{p} := p(x, y, z, t)(\mathbf{i}+\mathbf{j}+\mathbf{k})$. Thus we obtained the following equation

$$\frac{\partial}{\partial t}(\nabla \cdot \boldsymbol{\varphi}) + \nabla \cdot (\mathbf{u}\nabla \cdot \boldsymbol{\varphi}) = -\nabla \cdot \mathbf{p} + \frac{1}{Re}\nabla^2(\nabla \cdot \boldsymbol{\varphi}) \tag{12}$$

The integration of the above equation is in the form of

$$\frac{\partial \boldsymbol{\varphi}}{\partial t} + \mathbf{u}\nabla \cdot \boldsymbol{\varphi} = -\mathbf{p} + \frac{1}{Re}\nabla(\nabla \cdot \boldsymbol{\varphi}) + \mathbf{C} \tag{13}$$

here the integration parameter satisfies $\nabla \cdot \mathbf{C} = 0$. By noting the fact that $\nabla \cdot \mathbf{C} = \nabla \cdot \mathbf{u} = 0$, we may incorporate the integration parameter into the velocity term $\mathbf{u} \to \mathbf{u} + \mathbf{C}$ without changing the form of the equation. Write the above equation in the scalar forms

$$\frac{\partial \phi}{\partial t} + u_1 \nabla \cdot \boldsymbol{\varphi} = -p + \frac{1}{\text{Re}} \frac{\partial}{\partial x}(\nabla \cdot \boldsymbol{\varphi})$$
$$\frac{\partial \phi}{\partial t} + u_2 \nabla \cdot \boldsymbol{\varphi} = -p + \frac{1}{\text{Re}} \frac{\partial}{\partial y}(\nabla \cdot \boldsymbol{\varphi}) \quad . \tag{14}$$
$$\frac{\partial \phi}{\partial t} + u_3 \nabla \cdot \boldsymbol{\varphi} = -p + \frac{1}{\text{Re}} \frac{\partial}{\partial z}(\nabla \cdot \boldsymbol{\varphi})$$

Note that for incompressible condition, we also have $\nabla^2 \phi = 0$. By the summation of the above equations we finally obtain

$$\frac{\partial \phi}{\partial t} + \frac{1}{3}(\frac{\partial \phi}{\partial x} + \frac{\partial \phi}{\partial y} + \frac{\partial \phi}{\partial z})^2 - \frac{1}{\text{Re}}(\frac{\partial^2 \phi}{\partial x \partial y} + \frac{\partial^2 \phi}{\partial y \partial z} + \frac{\partial^2 \phi}{\partial z \partial x}) + p = 0 \qquad \blacksquare$$

In this way, we have reduced the NSE into the above simple scalar equation. Eq.(6) gives us a simpler expression for the dynamics of the potential field of fluids. The benefit of using such potential expression is that all the velocity components can be later on recovered by the derivation of the potential function. This equation indicates that the vector potential has no impact on the dynamics of the averaged potential field due to the vanishing of such component in the equation. In other words, the dynamics of the velocity potential is independent of its vector potential part. However, the total vector field can be recovered later on by superimposing the velocity potential and vector potential parts. In the next section, we are going to apply Lie group analysis on the equation.

## Symmetry Analysis

Since looking at the domain and codomain of a given differential equation as a product manifold can be naturally generalized to jet manifolds, it provides a more unified and concise language to the Lie symmetry method by studying jet bundles and corresponding Cartan distributions of the differential equations, we plan to discuss the problem here using the machinery of jet bundle. We shall first introduce some basic definitions related in the present paper, readers who wish to explore more details can go to [18, 19]

Let $(E, \pi, M)$ be a locally trivial smooth bundle over a smooth manifold $M$, of which the set of *sections* is $\Gamma(\pi) = \{\sigma : \pi \circ \sigma = id_M\}$. The *r-jet bundle manifold* associated with the bundle can be defined as

$$J_p^r(\pi) = \{j_p^r \sigma : p \in M, \sigma \in \Gamma(\pi)\}$$

where the *r-th* equivalence class $j_p^r \sigma$ is called *r-jet* of $\sigma$ at $p$. Thus, the disjoint union of the fiber manifolds

$$J^r(\pi) = \bigsqcup_{p \in M} j_p^r(\pi)$$

forms a smooth *vector bundle* $(J^r(\pi), \pi_r, M)$ that is called *jet bundle*. The construction also induces *affine bundles* $(J^r(\pi), \pi_{r,k}, J^k(\pi))$, here $0 \leq k < r$ and $J^0(\pi) = E$.

Suppose $\theta \in J^r(\pi)$, Let $\Gamma_\sigma^r \subset J^r(\pi)$ be the graph of r-jets, then the span of all planes tangent to the graph $\Gamma_\sigma^r |_\theta$ is called the *Cartan plane* $\mathcal{C}_\theta$. The disjoint union of the Cartan planes

$$\mathcal{C} = \bigsqcup_{\theta \in J^r(\pi)} \mathcal{C}_\theta$$

is an integrable distribution that is called the *Cartan distribution*, which is the basic geometric structure on the manifold $J^r(\pi)$.

Given an partial differential equation as

$$F\left(t, x, y, z, \phi, \phi_x, \phi_y, \phi_z, \phi_{xy}, \phi_{yz}\right) = 0 .$$

Then the set

$$\mathcal{E} = \{(x^i, \phi_j, \phi_{ji}) : F\left(x^i, \phi_j, \phi_{ji}\right) = 0\}$$

defines a submanifold in $J^2(\pi) = (\mathbb{R} \times \mathbb{R}^3, \pi, \mathbb{R})$. The Cartan distribution on the 2-jet manifold in the local coordinates then can be characterized by the following contact forms

$$\omega_i = d\phi_i - \phi_{i+1} dx^i$$

of which the restriction on the equation $\mathcal{E}$ induces the Cartan distribution on the equation [18, 20]

$$\mathcal{C}_\theta(\mathcal{E}) = \mathcal{C}_\theta \cap T_\theta \mathcal{E}$$

where $T_\theta \mathcal{E}$ is the tangent space on the equation at $\theta$. A maximal integral manifold of the Cartan distribution $\mathcal{C}_\theta(\mathcal{E})$ is called the *general solution* to the equation $\mathcal{E}$.

As for the equation in this paper, we can define a nowhere vanishing vector field on $J^0(\pi)$ spanned by

$$X = \xi_1 \frac{\partial}{\partial x} + \xi_2 \frac{\partial}{\partial y} + \xi_3 \frac{\partial}{\partial z} + \xi_4 \frac{\partial}{\partial t} + \eta \frac{\partial}{\partial \phi}$$

which is also generated by a one-parameter group $\Phi \in C^\infty(J^0(\pi))$. Then $X$ is called the *Lie point symmetry* if its lifting on $J^2(\pi)$ is in the form of

$$X^{(2)} = X + \eta^i \frac{\partial}{\partial \phi_i} + \eta^{ji} \frac{\partial}{\partial \phi_{ji}}$$

such that

$$X^{(2)} F\left(x^i, \phi_j, \phi_{ji}\right)\big|_{F=0} = 0$$

where

$$\eta^i = D_i \eta - \phi^j D_i \xi_j$$
$$\eta^{ji} = D_j \eta^i - \phi^{jk} D_j \xi_k$$

and the operator of total differentiation is in the form of

$$D_i = \frac{\partial}{\partial x^i} + \phi_i \frac{\partial}{\partial \phi_i} + \phi_{ji} \frac{\partial}{\partial \phi_{ji}}.$$

Furthermore, if $\Phi \in C^\infty(J^1(\pi))$, we call $X$ the *Lie contact symmetry*. One can easily check that the lifting of Lie symmetries preserve the Cartan distribution [18, 20], e.g.

$$\mathcal{L}_{X^{(2)}} \omega = \lambda \omega, \qquad \forall \lambda \in C^\infty\left(J^1(M)\right)$$

Recall the one-parameter Lie transformation group of the form [21]

$$\begin{aligned}
x^* &= x + \varepsilon \xi_1(x, y, z, t, \phi) + O(\varepsilon^2) \\
y^* &= y + \varepsilon \xi_2(x, y, z, t, \phi) + O(\varepsilon^2) \\
z^* &= z + \varepsilon \xi_3(x, y, z, t, \phi) + O(\varepsilon^2) \\
t^* &= t + \varepsilon \xi_4(x, y, z, t, \phi) + O(\varepsilon^2) \\
\phi^* &= \phi + \varepsilon \eta_1(x, y, z, t, \phi) + O(\varepsilon^2)
\end{aligned} \qquad (15)$$

and the corresponding infinitesimal generator

$$X = \xi_1 \frac{\partial}{\partial x} + \xi_2 \frac{\partial}{\partial y} + \xi_3 \frac{\partial}{\partial z} + \xi_4 \frac{\partial}{\partial t} + \eta \frac{\partial}{\partial \phi} \tag{16}$$

where

$$\xi_i = \xi_i(x^i, \phi) = \frac{dx^{i*}}{d\varepsilon}\Big|_{\varepsilon=0}$$
$$\eta = \eta(x^i, \phi) = \frac{d\phi^*}{d\varepsilon}\Big|_{\varepsilon=0} \tag{17}$$

and denote $(x^i) = (x, y, z, t)$, $i = 1, 2, 3, 4$. Also, the second prolongation (lifting) of the vector field $X$ can be written as

$$X^{(2)} = X + \eta^i \frac{\partial}{\partial \phi_i} + \eta^{ji} \frac{\partial}{\partial \phi_{ji}} \tag{18}$$

By applying

$$X^{(2)} F(x^i, \phi_j, \phi_{ji})\big|_{F=0} = 0 \tag{19}$$

the determining equations for $\xi_i$ and $\eta$ can be obtained. Usually, solving the determining equations is not easier than solving the original differential equations. It was the reason why the Lie group method was not of much practical use before the developing of modern computers. Right now by using CAS software or related packages, we can readily solve the systems of linear partial differential equations. For the scalar equation Eq.(6), the determining equations are obtained and solved by using MAPLE.

## Results

The corresponding Lie subalgebras and invariant solutions are presented case by case as follows.

**Case 1:** $p = 0$. In Eq.(6), we can treat the pressure as the source term with setting the pressure to be zero for the simplest case.

$$\xi_1 = \theta x + \zeta t + \kappa(y - z) + \alpha$$
$$\xi_2 = \theta y + \zeta t + \kappa(z - x) + \beta$$
$$\xi_3 = \theta z + \zeta t + \kappa(x - y) + \gamma$$
$$\xi_4 = 2\theta t + \delta$$
$$\eta_1 = \frac{3}{2}\zeta x$$

where the Greek alphabets are arbitrary constants. By setting each constant coefficient to be zero sequentially, the Lie subalgebras of the equation, in this case, are spanned by

$$X_1 = \frac{\partial}{\partial x}, \ X_2 = \frac{\partial}{\partial y}, \ X_3 = \frac{\partial}{\partial z}, \ X_4 = \frac{\partial}{\partial t}$$

$$X_5 = t\frac{\partial}{\partial x} + t\frac{\partial}{\partial y} + t\frac{\partial}{\partial z} + \frac{3x}{2}\frac{\partial}{\partial \phi}$$

$$X_6 = (y-z)\frac{\partial}{\partial x} + (z-x)\frac{\partial}{\partial y} + (x-y)\frac{\partial}{\partial z}$$

$$X_7 = t\frac{\partial}{\partial t} + \frac{x}{2}\frac{\partial}{\partial x} + \frac{y}{2}\frac{\partial}{\partial y} + \frac{z}{2}\frac{\partial}{\partial z}$$

and the possible invariant solutions

$$\phi = A, \ \phi = \frac{3x^2}{4t} + A, \ \phi = \frac{3y^2}{4t} + A, \ \phi = \frac{3z^2}{4t} + A$$

$$\phi = \frac{3z^2}{4t} + B\mathrm{erf}\left(\frac{\sqrt{-\mathrm{Re}}|z-y|}{2\sqrt{t}}\right) + A$$

$$\phi = \frac{3z^2}{4t} + B\mathrm{erf}\left(\frac{\sqrt{-\mathrm{Re}}|z-x|}{2\sqrt{t}}\right) + A$$

$$\phi = B + \int_{z/y} \frac{3da}{2\mathrm{Re}a\ln(a) - \mathrm{Re}a^2 - 3Aa + \mathrm{Re}}$$

.

**Case 2:** $p = p_0 = \mathrm{const}$. Using the same procedure (hereafter we omit the presentation of infinitesimal generators) we obtain the following Lie subalgebras

$$X_1 = \frac{\partial}{\partial x}, \ X_2 = \frac{\partial}{\partial y}, \ X_3 = \frac{\partial}{\partial z}, \ X_4 = \frac{\partial}{\partial t}$$

$$X_5 = t\frac{\partial}{\partial x} + t\frac{\partial}{\partial y} + t\frac{\partial}{\partial z} + \frac{3x}{2}\frac{\partial}{\partial \phi}$$

$$X_6 = (z-x)\frac{\partial}{\partial \phi}, \ X_7 = (y-x)\frac{\partial}{\partial \phi}$$

$$X_8 = (z-y)\frac{\partial}{\partial x} + (x-z)\frac{\partial}{\partial y} + (y-z)\frac{\partial}{\partial z}$$

$$X_9 = t\frac{\partial}{\partial t} + \frac{x}{2}\frac{\partial}{\partial x} + \frac{y}{2}\frac{\partial}{\partial y} + \frac{z}{2}\frac{\partial}{\partial z} - \frac{p_0 Re(x-y)^2}{2}\frac{\partial}{\partial \phi}$$

and the possible invariant solutions

$$\phi = A - p_0 t, \ \phi = A - x\sqrt{-3p_0}$$
$$\phi = A - y\sqrt{-3p_0}, \ \phi = A - z\sqrt{-3p_0}$$
$$\phi = \frac{3z^2}{4t} - p_0 t + A$$

.

**Case 3:** $p = p(t)$. In this case, we assume that the pressure depends on time only. The Lie subalgebras of the equation for this case are

$$X_1 = \frac{\partial}{\partial x}, \ X_2 = \frac{\partial}{\partial y}, \ X_3 = \frac{\partial}{\partial z}, \ X_4 = \frac{\partial}{\partial t} - p(t)\frac{\partial}{\partial \phi}$$

$$X_5 = t\frac{\partial}{\partial x} + t\frac{\partial}{\partial y} + t\frac{\partial}{\partial z} + \frac{3x}{2}\frac{\partial}{\partial \phi}$$

$$X_6 = (z-x)\frac{\partial}{\partial \phi}, \ X_7 = (y-x)\frac{\partial}{\partial \phi}, \ X_8 = \frac{\partial}{\partial \phi}$$

$$X_9 = (z-y)\frac{\partial}{\partial x} + (x-z)\frac{\partial}{\partial y} + (y-z)\frac{\partial}{\partial z}$$

$$X_{10} = t\frac{\partial}{\partial t} + \frac{x}{2}\frac{\partial}{\partial x} + \frac{y}{2}\frac{\partial}{\partial y} + \frac{z}{2}\frac{\partial}{\partial z} - tp(t)\frac{\partial}{\partial \phi}$$

and the possible invariant solutions

$$\phi = A - \int p(t)dt$$

$$\phi = A - \int p(t)dt + \frac{3z^2}{4t}.$$

**Case 4:** $p = p(x)$. In the last, let's consider the situation that the pressure depends $x$ only, three translation groups are corresponding to the following Lie subalgebras

$$X_1 = \frac{\partial}{\partial t}, \quad X_2 = \frac{\partial}{\partial y}, \quad X_3 = \frac{\partial}{\partial z}$$

and the possible invariant solutions are in the form of

$$\phi = A - \sqrt{3}\int\sqrt{-p(x)}dx$$
$$\phi = A + \sqrt{3}\int\sqrt{-p(x)}dx.$$

## Discussion

From the above Lie group analysis, we have obtained the Lie subalgebras for some special cases, as well as their invariant solution. Although the physical meaning of each the solution has not yet been fully exploited, the approach in this paper has shed light on how to simplify coupled differential equations. By the summation of each component in the equations, we keep only the scalar information, which is more important in most of the cases, e.g. in wave equation the amplitude of the wave and in heat equation the intensity of the heat.

If we specify the pressure $p(t) = (x + y + z)t$ we can construct the invariant solutions for this special case as

$$\phi = A - \frac{t^2}{2}(x + y + z) - \frac{3t^5}{20}$$

$$\phi = \frac{3z^2}{4t} - \frac{t^2}{2}(y + x) - \frac{t^5}{15} + A$$

as well as $p(t) = \sin(xyz - t)$, and the corresponding invariant solutions

$$\phi = A - \cos(t - xyz)$$

$$\phi = A + \frac{3z^2}{4t} - \cos(t - xyz)$$

.

With the aid of the Stokes-Helmholtz decomposition, we have uncoupled the 3-dimensional Navier-Stokes equation into a scalar equation of the velocity potential. This equation gives a neat way to describe the dynamical behaviors of flows in the toroidal field. The approach of simplification indicates an averaging magnitude of the velocity vectors. The vector potential does not affect the dynamics of the average potential field since it is vanished in the equation after averaging. The Lie group analysis has helped us to construct some invariant solutions to the equation within some instances.

We applied an averaged velocity potential to simplify the NSE in the paper, of which the technique may also help simplify other similar partial differential equations.